\documentclass{article}
\usepackage{latexsym}

\usepackage{amsmath,amssymb,amsfonts,array}

\def\CC{\mathbb{C}}

\def\RR{\mathbb{R}}
\def\ZZ{\mathbb{Z}}

\author{Klaus Mohnke\thanks{Universit\"at-Gesamthochschule Siegen,
e-mail:
mohnke@mathematik.uni-siegen.de} }
\title{Lagrangian Embeddings in the
Complement of Symplectic Hypersurfaces}
\date{}
\newtheorem{theorem}{Theorem}[section]
\newtheorem{definition}[theorem]{Definition}
\newtheorem{proposition}[theorem]{Proposition}
\newtheorem{lemma}[theorem]{Lemma}
\newtheorem{remark}[theorem]{Remark}
\newtheorem{corollary}[theorem]{Corollary}
\begin{document}
\bibliographystyle{plain}

\maketitle

\begin{abstract}
For  a given embedded Lagrangian in the complement of a complex
hypersurface we show existence of a holomorphic disc in the
complement having boundary on that Lagrangian.
\end{abstract}

\section{Introduction and Main Results}

Let $L\subset\CC^n$ be a Lagrangian, i.e.~immersed (real)
$n$--manifolds $L\subset\CC^n$ such that the restriction
$\omega|_L=0$ for the symplectic $2$-form $\omega = d\theta =
\sum_{i=1}^ndx^i\wedge dy^i$ vanishes, where $\theta =
\sum_{i=1}^nx^idy^i$ is the Liouville form and $z^i\!=\!
x^i\!+\!\mbox{i}y^i$ are the (complex) coordinates on $\CC^n$. Hence
$\theta|_L$ is a closed $1$-form and thus defines an element in
${\bf H^1}(L,\ZZ)$. In \cite{Gromov} Gromov constructed non
constant holomorphic discs with boundary on any {\em embedded},
closed Lagrangian in $\CC^n$. He concluded that such a Lagrangian
is never exact, i.e. $[\theta|_L]\in {\bf H^1}(L;\ZZ)$ defines a
non vanishing class in de Rham cohomology: the symplectic volume of
the disc is positive by the compatibility of complex and symplectic
structure and so is the integral of $\theta$ over the boundary of
the disc by Stokes' theorem. In this paper we specialize the disc
found by Gromov in the presence of a complex hypersurface. Let us
first fix some notation (see
e.g.~\cite{Audin/Lalonde/Polterovich}).

\begin{definition}
A complete Riemannian manifold is called {\em geometrically bounded
(g.b.)}~if its sectional curvature is bounded above and the
injectivity radius is bounded below. A submanifold is  g.b.~if its
second fundamental form is bounded above and if there are constants
$\epsilon,C>0$ such that if two points of the submanifold can be
joined by a geodesic of length $l\le \epsilon$ in the ambient space
they can be joined by a geodesic of length $Cl$ in the submanifold.
\end{definition}

Note that if the manifold (submanifold) in question is compact it is,
of course, geometrically bounded. Now we are ready to state the main
statements of this note.

\begin{theorem}\label{lagrange}
(i) Suppose $L\cong\RR\times{\bf S^1}$ is a g.b.~Lagrangian
embedding into $\CC^2$ such that there exists a complex line $H$ in
its complement nontrivially linked by $L$,
i.e.~$0\neq\mbox{lk}(H,.):{\bf H_1}(L)\to\ZZ$.
Then each complex linear function which is
non-constant on $H$ is unbounded on $L$.\\ (ii) There is no
Lagrangian embedding of the Klein bottle into the complement of a
(real) plane in $\CC^2$, nontrivially linking it.
\end{theorem}

\begin{remark}
(0) The Klein bottle is the only closed surface for which it is
unknown wether it admits a Lagrangian embedding into $\CC^2$ or
not.  For this some progress has been reported by Hofer and
Luttinger \cite{Hofer:Report}.\\ (1) We may construct embedded
totally real Klein bottles in the complement of a plane
nontrivially linking it. Also note that there is an immersed
Lagrangian Klein bottle in the complement of a plane linking it
nontrivially (see Exercise~1.4.2.~of
\cite{Audin/Lalonde/Polterovich})\\ (2) Lalonde proves in
\cite{Lalonde} that if there exists an embedded Lagrangian Klein
bottle in $\CC^2$ then the double points of the immersed one in (1)
can be suppressed by a symplectic Whitney disc. It follows from
Theorem~\ref{lagrange} that this disc may never lie in the
complement of the plane.
\end{remark}

All this will follow from the main result of the paper. We will
formulate it in a more general setting refining the arguments to
specify the holomorphic disc to lie in the complement of a complex
hypersurface with boundary on a given Lagrangian in that
complement.


\begin{theorem}\label{main}
Let $M^{2n}$ be a symplectic manifold and $\Sigma\subset M$ be a
symplectic hypersurface such that there exists a compatible  almost complex
structure such that $M$ and $\Sigma$ are geometrically bounded and
there is no nonconstant
holomorphic sphere in $M$ with nonpositive intersection index with
the hypersurface $\Sigma$ inside $M$. Suppose $U\subset M$ is a
neighbourhood of $\Sigma$ and $L\subset (M\setminus U)\times\CC$ is
a g.b.~Lagrangian such that its projection onto the second
component is bounded. Then there exists a loop $\gamma\subset L$
which bounds a nonconstant holomorphic disc in
$(M\setminus\Sigma)\times \CC$.
\end{theorem}

{\em Proof of Theorem~\ref{lagrange}.} (i) Assume there is a
bounded complex linear function which is nonconstant on $H$ then
the hypotheses of Theorem~\ref{main} are satisfied and in the
complement of $H$ there is a disc with boundary on $L$ and
nonvanishing symplectic area. However, this contradicts
$\mbox{lk}(H,.)\neq 0$ since the boundary of the disc represents an
element of its kernel which is not torsion.\\ (ii)  Let $L$ be a
Lagrangian embedding of a closed surface into the complement of a
(real) plane. 
Since for a plane to be symplectic is an open condition
we may assume that the plane is symplectic (by
replacing it with a nearby plane) and fix a constant, and thus
integrable complex structure on $\CC^2$ such that the plane is a
complex line $H\subset\CC^2$. In the complement of $H$ we find a
disc with boundary on $L$ and nonvanishing symplectic area. Thus we
obtain an element in the kernel of $\mbox{lk}(H,.):{\bf
H_1}(L)\longrightarrow{\bf Z}$ which is not torsion. Hence
$b_1(L)\ge 2$ and the statement follows. $\Box$

There are also consequences in higher dimensions which are
interesting from another point of view: Here we have much more
space to construct Lagrangians of all kinds, for e.g. via a
''cabeling'' procedure (see Proposition~1.2.3 in
\cite{Audin/Lalonde/Polterovich}). But the following corollary
shows that not everything is possible:
\begin{corollary}
Let $L$ be a g.b.~Lagrangian in the complement of a complex
hyperplane $H\subset\CC^n$. Then either\\ (i) each complex linear
function on $\CC^n$ is unbounded on $L$ or\\ (ii) the topology of
$L$ is restricted by $b_1(L)\ge 1$\\
\hspace*{0.5cm} $\bullet$ if $L$ and the plane are unlinked or\\
\hspace*{0.5cm} $\bullet$ if each complex linear function which
is bounded on $L$ vanishes on $H$\\ (iii) or $b_2(L)\ge 2$.
\end{corollary}

{\em Proof.} The hypotheses of Theorem~\ref{main} are satisfied.
This, together with the assumptions on the linking and the
existence of the complex linear function, implies the statement as
in the proof of Theorem~\ref{lagrange}. If $L$ is closed we may
take any real $4$-plane by the perturbation argument of that proof.
$\Box$

\begin{corollary}
In particular, if $L$ is closed and nontrivially linked with any real
$2(n\!-\!1)$-plane, then $b_1(L)\ge 2$. That means, for $n\ge 3$,
that a Lagrangian embedding of ${\bf S^{n-1}}\times{\bf S^1}$ in
$\CC^n$ (which e.g.~can be constructed from Whitney's Lagrangian
immersion of the sphere ${\bf S^{n-1}}$ into $\CC^{n-1}$ via the
cabling mentioned above) trivially links any real $2(n-1)$-plane
in its complement. $\Box$
\end{corollary}

\section{Proof of Theorem~\ref{main}}\label{poltero}

We prove the observation refining the arguments used to prove the
corresponding result of Gromov \cite{Gromov} (for more details see
e.g.~\cite{Audin/Lalonde/Polterovich}).

Denote by $\pi:M\times\CC\longrightarrow M$ the projection onto the
first component. Fix an $\alpha>2$. We define

\begin{eqnarray*}
{\cal F} & := & \{u\in{\bf C}^\alpha(({\bf
\Delta},\partial{\bf\Delta},1),(M\times\CC,L,p))\mid
u\cong p \mbox{ in } \pi_2(M\times\CC,L,p)\}\\ {\cal G} & := &
\{g\in{\bf C}^{\alpha-1}({\bf \Delta}\times
(M\times\CC), T^{\CC}(M\times\CC)\mid
\pi_\ast(g|_{{\bf\Delta}\times U\times\CC})=0\}\\
{\cal M} & := &
\{(u,g)\in {\cal F}\times{\cal G}\mid \overline{\partial}u-g\},
\end{eqnarray*}
where by $\overline{\partial}$ denotes the Cauchy--Riemann operator
$\frac{1}{2}(\frac{\partial}{\partial x}+J\frac{\partial}{\partial y})$
in complex coordinates $z=x+\mbox{i}y$ on the disc.

The linearization
$$
D_{(u,g)}: {\bf C}^\alpha(({\bf
\Delta},\partial{\bf\Delta},1),(u^\ast TM\times\CC,u^\ast TL,0))
\times {\cal C}^{\alpha-1}({\bf\Delta},u^\ast TM\times \CC)
\longrightarrow {\cal C}^{\alpha-1}({\bf\Delta},u^\ast TM\times \CC)
$$
of the equation $\overline{\partial}u=g$ is surjective for the
following reason: there are sufficient degrees of freedom coming
from the second component to perturb the equation at the boundary
of the disc. Hence an element in the cokernel would satisfy the
Cauchy-Riemann equation with Dirichlet conditions on the boundary
which imply its vanishing. Thus

\begin{proposition}\label{smale}
${\cal M}\subset{\cal F}\times{\cal G}$ is a Banach submanifold.
The natural projection $P:{\cal M}\longrightarrow{\cal G}$ is a
Fredholm map of index $0$. It is regular at $(p,0)$. $\Box$
\end{proposition}

Note that $P^{-1}(0)=\{(p,0)\}$ consists of a single point. Pick
$h:=(0,B)$ according to the splitting
$T(M\times\CC)\cong \pi^\ast TM\times\underline{\CC}$, $\underline{\CC}$
being the trivial complex line bundle, with
$B\in\CC$ understood as a constant section of it.
Connect $h$ it to $0$ by a path $\gamma\subset{\cal
G}$ transversal w.r.t.~$P$. Hence $P^{-1}(\gamma)$ is a
one-dimensional manifold with boundary consisting of $(p,0)$ and
$P^{-1}(g)$. Via Cauchy-integral or Poisson-formula one proves that
there is no solution to $\overline{\partial}u = h$ provide $B$ is
sufficiently large. (see e.g.~\cite{Audin/Lalonde/Polterovich}).
Therefore $P^{-1}(g)=\emptyset$. Thus $P^{-1}(\gamma)$ has to be
non-compact.

Let $t_0 := \inf\{t\mid P^{-1}(\gamma(t))=\emptyset\}$.  Pick a
sequence $t_n\uparrow t_0>0$ and $(u_n,g_n)\in
P^{-1}(\gamma(t_n))\subset{\cal M}$ with $\lim_{n\to\infty}g_n\in
g\in{\cal G}$ which does not converge in ${\cal M}$. From
\cite{Gromov} (for a longer discussion see also
\cite{Audin/Lalonde/Polterovich})
we therefor know that there is ''bubbling'': a
subsequence converges in the sense of Gromov to a solution of
$\overline{\partial}u=g$ with
$u:({\bf\Delta},\partial{\bf\Delta})\longrightarrow
(M\times\CC,\Sigma\times\CC)$ representing a different homotopy
class then the $u$'s considered so far and a (nonempty!) collection
of nonconstant holomorphic discs with boundary on $L$ and
nonconstant holomorphic spheres counted with multiplicity. The convergence
is such that the
homology class of the $u_n$ is preserved. Hence the sum of intersection numbers
of these holomorphic curves and $u$ with $\Sigma\times\CC$ is zero. But the discs
have non-negative intersection number with the hypersurface for
they are holomorphic and the boundary is disjoint from it due to
Lemma~\ref{intersection}. 
Also
the spheres have positive intersection number with $\Sigma\times\CC$
by assumption.

Thus it suffices to show that $u$ has  nonnegative intersection
number with $\Sigma\times\CC$. For this excludes spheres in the
collection of holomorphic curves. Hence there {\em is} a nonconstant
holomorphic disc with boundary on $L$ and zero intersection index
with $\Sigma\times\CC$ which implies that it has to lie in the
complement.

Let us prove the statement about $u$. Assume the contrary. Pick a point
$q\in u({\bf
\Delta})\cap \Sigma\times\CC$. We claim that its intersection index is
positive. Since $\pi_\ast(g|_{{\bf \Delta}\times U\times\CC})=0$
we obtain a holomorphic map $\pi\circ u:
u^{-1}(U\times\CC)\subset{\bf\Delta}\to M$
near the preimage of $q$. The intersection indices of $u$ with
$\Sigma\times\CC$ at $q$ and
$\pi\circ u$ with $\Sigma$ at $\pi(q)$ agree. Due to Lemma~\ref{intersection}
the latter is a
well--defined (finite!) positive number unless the whole image of
that neighbourhood lies inside $\Sigma$. This cannot happen because
then $u$ would map completely into $\Sigma\times\CC$. But at least
the boundary lies in $L\subset (M\setminus \Sigma)\times\CC$.
$\Box$

\begin{remark}
Theorem~\ref{lagrange} has a simpler proof directly applying Gromov's original theorem
using the fact that
$\CC^n\setminus\CC^{n-1}$ can be symplectomorhically embedded in $(\RR\times{\bf
S^1})\times\CC^{n-1}$

This situation belongs to the case when the complement of $\Sigma$ is Stein. Here
one may argue with the non-existence of holomorphic curves touching some
(weakly) pseudoconvex boundary from the inside instead. However, the statement
is true in a much more general context.
\end{remark}

Let us finish with a discussion on the positiveness of the intersection
indices of an embedded almost complex hypersurface and a (pseudo)holomorphic
curve not lying inside the former.
\begin{lemma}\label{intersection}
Consider an almost comlex manifold $M$.
Let $\Sigma\subset M$ be an embedded almost complex hypersurface and 
$u:{\bf\Delta}\longrightarrow M$
be a (pseudo)holomorphic curve. Then
either the image of $u$ is completely contained inside $\Sigma$ or
$u^{-1}(\Sigma)\subset{\bf\Delta}$ is a discrete set inside the
open disc ${\bf\Delta}$ each point coming with a finite intersection
index.
\end{lemma}
{\em Proof.} The notion of the index will become clear during the discussion
of the proof. The following  argument is due to Brian White.
Near an intersection point we may pick coordinates such that $M$ looks
like a subset of $\CC^{n-1}\times\CC$ (containing $0$ say) and the
almost complex structur has the form $J(z,\zeta)=\left(\begin{array}{cc} j(z,\zeta) & 0\\
O(\zeta)(z,\zeta) & i \end{array}\right)$ where $j$ is an almost complex structure
on $\CC^{n-1}$, $i$ is the ordinary complex multiplication and $O(\zeta)$
is a smooth $2(n-1)\times 2$-matrix-valued function (with real entries)
vanishing of first order with $|\zeta|$ satisfying $O(\zeta)j+iO(\zeta)=0$.
Now, if $(U,u):{\bf\Delta}\longrightarrow \CC^{n-1}\times\CC$ is a pseudo-holomorphic curve through
$(0,0)$ then in the coordinate chart it satisfies the following differential
equation
$$
\frac{\partial u}{\partial \overline{z}}=O'(u)DU,
$$
again with a matrix-valued function $O'(\zeta)$ of the kind above. Differentiating this
equation we end up with the differential inequality
$$
|{\bf \Delta}u|\le K(|u|+|Du|)
$$
for the $u$-part of that pseudoholomorphic curve. As in \cite{Micallef/White} 
we  conclude that either the whole image of the curve will be contained in
the hypersurface or there is a positive integer $m$ such that $u(\xi)=\xi^m+O(|\xi|^{m+1})$
proving the statement of the lemma. Another argument was given by Jean--Claude Sikorav.
He concludes the result using similarity principle at the first order differential equation for $u$
as in his paper \cite{Sikorav}.
$\Box$

{\em Acknowledgements.} I would like to thank Yakov Eliashberg for
the initial discussions on this subject, Leonid Polterovich for his
cordial and patient help with the paper, Brian White and Jean-Claude Sikorav
for the help with the proof of Lemma~\ref{intersection}.

\end{document}